\documentclass[aos,MSNbibl,nameyear,rotating,dvips]{arximspdf}
\usepackage{dcolumn}

\doi{10.1214/14-AOS1225} 
\volume{42}
\issue{4}
\pubyear{2014}
\firstpage{1347}
\lastpage{1360}

\makeatletter
\newcommand{\OA}{\operatorname{OA}}
\newcommand{\SOA}{\operatorname{SOA}}
\newcommand{\GOA}{\operatorname{GOA}}
\newcolumntype{d}[1]{D{.}{.}{#1}}
\newtheorem{lem}{Lemma}
\newtheorem{teo}{Theorem}
\newtheorem{prop}{Proposition}
\newproclaim{defin}{Definition}
\newcommand{\comb}[2]{\pmatrix{#1\cr #2}}
\newcommand{\tcomb}[2]{{#1\choose #2}}
\makeatother

\begin{document}
\begin{frontmatter}

\title{A characterization of strong orthogonal arrays of~strength three\thanksref{T1}}
\runtitle{Strong orthogonal arrays}

\begin{aug}
\author[a]{\fnms{Yuanzhen}~\snm{He}\ead[label=e1]{heyuanzhen@gmail.com}}
\and
\author[b]{\fnms{Boxin}~\snm{Tang}\corref{}\ead[label=e2]{boxint@sfu.ca}}
\runauthor{Y. He and B. Tang}
\affiliation{Chinese Academy of Sciences and Simon Fraser University}
\address[a]{Academy of Mathematics\\
\quad and Systems Science\\
Chinese Academy of Sciences\\
Beijing 100190\\
China\\
\printead{e1}}
\address[b]{Department of Statistics\\
\quad and Actuarial Science\\
Simon Fraser University\\
Burnaby, British Columbia V5A 1S6\\
Canada\\
\printead{e2}} 
\end{aug}
\thankstext{T1}{Supported by the Natural Sciences and Engineering Research Council of Canada.}

\received{\smonth{7} \syear{2013}}
\revised{\smonth{12} \syear{2013}}

%
\begin{abstract}
In an early paper, He and Tang [\textit{Biometrika} \textbf{100} (2013) 254--260]
introduced and studied a new
class of designs,
strong orthogonal arrays, for \mbox{computer} experiments, and characterized
such arrays through generalized orthogonal arrays. The current paper
presents a simple characterization for strong orthogonal arrays of
strength three.
Besides being simple, this new characterization through a notion of
semi-embeddability is more direct and penetrating in terms of revealing
the structure of strong orthogonal arrays. Some other results on strong
orthogonal arrays of strength three are also obtained along the way, and
in particular, two $\operatorname{SOA}(54, 5, 27, 3)$'s are constructed.
\end{abstract}

%
\end{frontmatter}

\section{Introduction}\label{sec1}

Computer models are powerful tools that enable researchers
to investigate complex systems from almost every imaginable field of
studies in natural sciences,
engineering, social sciences and humanities.
Computer models can be stochastic or deterministic; we consider
deterministic computer models. When the computer code representing
a computer model is expensive to run, it is desirable to build a cheaper
surrogate model. Computer experiments are concerned with
the building of a statistical surrogate model based on the data
consisting of
a set of carefully selected inputs and the corresponding outputs
from running a computer code.

Designing a computer experiment, that is, the selection of inputs, is
a crucial step in the process of model building.
No matter how elaborate and sophisticated a model building process is,
a statistical model contains no more information than what the data can offer.
In the past two decades, space-filling designs have been widely accepted
as appropriate designs for computer experiments. A space-filling design
refers to, in a very broad sense, any design that strews its points
in the design region in some uniform fashion. The uniformity of a design
may be evaluated using a distance criterion
[\citet{JohMooYlv90}] or a discrepancy criterion
[\citet{FanMuk00}]. See \citet{SanWilNot03}
and \citet{FanLiSud06} for more details.
Orthogonality has also been playing a
significant role in constructing designs for computer experiments
as it, in addition to being useful in its own right, provides a stepping
stone to achieving uniformity [\citet{LinMukTan09}].

The curse of dimensionality, however, makes it extremely difficult
for the points of a design to provide a good coverage of
a high dimensional design region. Even 10,000 points are not enough for
a $2^m$ grid in an $m=14$ dimensional space, not to mention that
such a regular grid leaves a deep hole in the center of the $2^m$ points.
In such situations, it makes more sense to consider designs that
are space-filling in lower dimensional projections of the input space.
The idea of Latin \mbox{hypercube} designs is to achieve the maximum uniformity
in all one-dimensional projections [\citet{McKBecCon79}].
$\OA$-based Latin hypercubes [\citet{Tan93}] carry this idea further,
which give designs that, in addition to being Latin hypercubes,
achieve uniformity in $t$-dimensional margins when orthogonal arrays
of strength $t$ are employed.
One could also use orthogonal arrays directly [\citet{Owe92}] but
such designs do not perform well in one-dimensional projections
when orthogonal arrays have small numbers of levels.

\citet{HeTan13} introduced, constructed and studied a new class of arrays,
\emph{strong orthogonal arrays}, for computer experiments.
A strong orthogonal array of strength $t$ does as well as
a comparable orthogonal array in $t$-dimensional projections,
but the former achieves uniformity on finer grids than the latter
in all $g$-dimensional projections for any $g \leq t-1$.
Consequently, Latin hypercubes constructed from a strong orthogonal array
of strength $t$ are more space-filling than comparable $\OA$-based Latin
hypercubes in all $g$-dimensional projections for any $2 \leq g \leq t-1$.
The concept of strong orthogonal arrays is motivated by the notion
of nets from quasi-Monte Carlo methods [\citet{Nie}].
The formulation of this new
concept has two advantages. First, strong orthogonal arrays are more
general than nets in terms of run sizes; and second, strong orthogonal
arrays are defined in the form and language that are familiar to
design practitioners and researchers. This not only makes
existing results from nets more accessible to design community
but also allows us to obtain new designs and theoretical
results.

The present article focuses on strong orthogonal arrays of strength three.
Through a notion of semi-embeddability, we provide a complete yet very
simple characterization for such arrays.
Though the characterization using generalized orthogonal arrays in \citet{HeTan13} is general, our new characterization for strength three
is more direct and revealing. Apart from this main result,
some other results on strong orthogonal arrays of strength three
are also obtained. In particular, we construct two
strong orthogonal arrays of 54 runs, five factors, 27 levels
and strength three.

The paper is organized as follows. Section~\ref{sec2} introduces some
notation and background material.
In Section~\ref{sec3}, a notion of semi-embeddability is defined,
through which we present the main result of the paper,
stating that a strong orthogonal array of
strength three exists if and only if a semi-embeddable orthogonal array of
strength three exists. We then examine the semi-embeddability and
nonsemi-embeddability of some orthogonal arrays. Section~\ref{sec4} constructs
two $\SOA(54,5,27,3)$'s. A discussion is given in Section~\ref{sec5}.

\section{Notation and background}\label{sec2}

This section provides a preparation for the rest of the paper
by introducing necessary notation and some background material.
An $n \times m$ matrix $A$ with its $j$th column taking levels
$0, 1, \ldots, s_j-1$ is said to be an orthogonal array of size $n$,
$m$ factors, and strength $t$ if for any $n \times t$ sub-matrix of~$A$,
all possible level combinations occur equally often.
Such an array is denoted by $\OA(n, m, s_1 \times\cdots\times s_m, t)$
in this paper. If at least two $s_j$'s are unequal,
the array is said to be asymmetrical or have mixed levels.
When $s_1 = \cdots= s_m =s$, we obtain a symmetrical orthogonal array, in
which case, the array is denoted by $\OA(n,m,s,t)$.
Since they were first introduced by \citet{Rad47}, orthogonal arrays have
been playing a prominent role in both statistical and combinatorial design
literature, and have become the backbone of designs for
multi-factor experiments. \citet{DeyMuk99} discussed
the construction and optimality of orthogonal
arrays as fractional factorial designs.
For a comprehensive treatment of orthogonal arrays, we refer to
\citet{HedSloStu99}.

Motivated by the notion of nets from quasi-Monte Carlo methods
[\citet{Nie}], \citet{HeTan13} introduced strong orthogonal
arrays.
Let $[x]$ denote the largest integer not exceeding $x$.
An $n \times m$ matrix with levels from
$\{0,1,\ldots,s^t-1\}$ is called a strong orthogonal array
of size $n$, $m$ factors, $s^t$ levels, and strength $t$ if any sub-array
of $g$ columns for any $g$ with $1 \leq g \leq t$ can be collapsed into an
$\OA( n,g, s^{u_1} \times s^{u_2} \times\cdots\times s^{u_g}, g)$
for any positive integers $u_1, \ldots, u_g$ with
$u_1 + \cdots+ u_g = t$, where collapsing into $s^{u_j}$ levels
is done by $[a/s^{t-u_j}]$ for $a=0, 1, \ldots, s^t-1$.
We use $\SOA(n,m,s^t,t)$ to denote such an array.
The following is an $\SOA(8,3,8,3)$:
\[
\left[\matrix{
0 & 0 & 0
\cr
2 & 3 & 6
\cr
3 & 6 & 2
\cr
1 & 5 & 4
\cr
6 & 2 & 3
\cr
4 & 1 & 5
\cr
5 & 4 & 1
\cr
7 & 7 & 7}\right],
\]
as we can easily check that:
\begin{longlist}[(iii)]
\item[(i)] The array becomes an $\OA(8,3,2,3)$ after the eight
levels are collapsed into two levels according to
$[a/4] = 0$ for $a=0,1,2,3$ and
$[a/4] = 1$ for $a=4,5,6,7$.

\item[(ii)] Any sub-array of two columns can be collapsed into
an $\OA(8,2, 2 \times4,2)$ as well as an $\OA(8,2, 4 \times2,2)$,
where collapsing into two levels is done by $[a/4]$ and
collapsing into four levels is done using $[a/2]$.

\item[(iii)] Any sub-array of one column is an $\OA(8,1,8,1)$.
\end{longlist}

\citet{Law96} introduced the concept of a generalized orthogonal
array. Extending a result of \citet{Law96},
\citet{HeTan13} showed that the existence of a strong orthogonal
array is equivalent to the existence of a generalized orthogonal array.
For the ease of presentation and the need of this paper,
here we give a \mbox{review} of this equivalence result only for the case
of strength three.
An $n \times(3m)$ matrix
$B= \{ (a_1, b_1, c_1); \ldots; (a_m, b_m, c_m) \}$
with entries from $\{0,1,\ldots,s-1\}$,
where, as indicated, the $3m$ columns are put into $m$ groups of
three columns each,
is called a generalized orthogonal array of size $n$, $m$ constraints,
$s$ levels and strength three if all the following matrices are
orthogonal arrays of strength three:
$(a_i, a_j, a_k)$ for any $1\leq i < j < k \leq m$,
$(a_i, b_i, a_j)$ for any $1\leq i \neq j \leq m$ and
$(a_i, b_i, c_i)$ for any $1\leq i \leq m$.
We use $\GOA(n,m,s,3)$ to denote such an array.

\begin{lem}\label{lemma1}
Let
$B= \{ (a_1, b_1, c_1); \ldots; (a_m, b_m, c_m) \}$
be a $\GOA(n,m,s,3)$. Define
%
\begin{equation}\label{eq1}
d_i = a_i s^2 + b_i s +
c_i.
\end{equation}
Then
$D= (d_1, \ldots, d_m)$ is
an $\SOA(n,m,s^3,3)$.
Conversely, if $D= (d_1, \ldots, d_m)$ is an $\SOA(n,m,s^{3\!},3)$,
then $B= \{ (a_1, b_1, c_1); \ldots; (a_m, b_m, c_m) \}$
is a $\GOA(n,m,\break  s,3)$, where $a_i, b_i, c_i$
are uniquely determined by $d_i$ as given in (\ref{eq1}).
\end{lem}

A bit explanation helps understand how $a_i, b_i, c_i$ are
obtained from $d_i$ in the second\vspace*{2pt} part of Lemma~\ref{lemma1}.
Every integer $0 \leq x \leq s^3-1$
can be uniquely written as $x= x_1 s^2 + x_2 s + x_3$
for some integers $x_1, x_2, x_3$ with $0 \leq x_j \leq s-1$. Applying this
fact to every component of $d_i$, we obtain
$d_i = a_i s^2 + b_i s + c_i$ for unique vectors
$a_i, b_i, c_i$, all with entries from $\{ 0, 1, \ldots, s-1 \}$.

Strong orthogonal arrays provide a new class of
suitable designs for computer experiments.
A strong orthogonal array of strength $t$
enjoys better space-filling properties than
a comparable orthogonal array in all dimensions lower than $t$ while
retaining the space-filling properties of the latter
in $t$ dimensions. Strong orthogonal arrays are more
general than nets in terms of run sizes. They
are defined in the form and language that are familiar to
design practitioners and researchers, and thus help to make
the existing results from nets more accessible to design community.
More importantly, this new formulation of the net idea in terms of
orthogonal arrays allows new designs and results to be found,
as has been shown in \citet{HeTan13} and will be further
demonstrated in the next two sections of the present
paper.\looseness=-1

The rest of the section discusses strong orthogonal arrays in the broad
context of quasi-Monte Carlo methods. To approximate an integral, Monte
Carlo methods evaluate the integrand at a set of points selected randomly,
whereas quasi-Monte Carlo methods do so at a set of points selected in a
deterministic \mbox{fashion}.
Specifically, to approximate $\int_{[0,1]^m} f(x) \,dx $,
quasi-Monte Carlo methods use $\sum_{i=1}^n f(x_i)/n$ where
$x_1, \ldots, x_n$ are a set of points in $[0,1]^m$ that are selected
deterministically and judiciously. The Koksma--Hlawka inequality
[\citet{Nie}, Theorem 2.11] states that
\[
\Biggl| \sum_{i=1}^n
f(x_i)/n - \int_{[0,1]^m} f(x) \,dx \Biggr| \leq
V(f) D^*_n (P),
\]
where $V(f)$ is the bounded variation of $f$ in the sense of
Hardy and Krause; $D^*_n(P)$ is the star discrepancy of the set $P$
of points $x_1, \ldots, x_n$, which is defined as~the maximum absolute difference
between the uniform distribution function
and the empirical distribution function based on the point set.
According to this result, the set of points
for quasi-Monte Carlo methods should therefore be chosen to have a small
star discrepancy.
When an infinite sequence of points
is considered, we use $D^*_n (S)$ to denote the star discrepancy given by
the first $n$ points of the sequence. The\vspace*{1pt} best general lower bounds
on $D^*_n(P)$ and $D^*_n (S)$ are those of \citet{Rot54} stating that
$D^*_n (P) \geq C_m   n^{-1} (\log n)^{(m-1)/2} $ for a point set
and
$D^*_n (S) \geq C_m   n^{-1} (\log n)^{m/2} $ for an infinite sequence,
where $C_m$ is a constant
independent of~$n$. But it is widely believed, though yet to be proved,
that
%
\begin{equation}\label{eq2}
D^*_n (P) \geq C_m n^{-1} (\log
n)^{m-1},\qquad D^*_n (S) \geq C_m
n^{-1} (\log n)^{m}.
\end{equation}
Halton sequences and corresponding Hammersley point sets attain the lower
bounds in (\ref{eq2}), but the implied constants $C_m$ grow superexponentially as
$m \rightarrow\infty$ [\citet{Nie}, Chapter~4].
What makes $(t,m,s)$-nets and $(t,s)$-sequences
attractive is that they have much smaller implied constants while
satisfying the lower bounds in (\ref{eq2}). Moreover, $(t,m,s)$-nets and
$(t,s)$-sequences contain an orthogonal array structure, which was
pointed out by \citet{Owe95}
and used by \citet{HaaQia10} to construct nested
space-filling designs for multi-fidelity computer experiments.

In what follows, we write $(w,k,m)$-nets
for $(t,m,s)$-nets and $(w,m)$-sequences for $(t,s)$-sequences
so as to be consistent in our notation for
this paper.
An elementary interval in base $s$ is an interval in
$[0,1]^m$ of form
\[
E = \prod_{j =1}^m \biggl[ \frac{c_j}{s^{d_j}},
\frac{c_j + 1}{s^{d_j}} \biggr),
\]
where nonnegative integers $c_j$ and $d_j$
satisfy $0 \leq c_j < s^{d_j}$.
For $0 \leq w \leq k$,
a $(w, k, m)$-net in base $s$ is a set of $s^k$ points in $[0,1]^m$
such that every elementary interval in base $s$ of volume $s^{w - k}$
contains exactly $s^w$ points.
Nets and related $(w,m)$-sequences were first defined by \citet{Sob67}
for base $s=2$ and later
by \citet{Nie87} for general base $s$.

A deeper connection of nets with orthogonal arrays was established
by \citet{Law96} and independently by \citet{MulSch96}.
These authors showed that a $(w,k,m)$-net
is equivalent to a generalized orthogonal
array. Inspired by this equivalence result, \citet{HeTan13} proposed
and studied strong orthogonal arrays for computer experiments.
Unlike generalized
orthogonal arrays, strong orthogonal arrays are in the ready-to-use format
and directly capture the space-filling properties of $(w,k,m)$-nets.
The following result is from \citet{HeTan13}.

\begin{lem}\label{lemma2}
If $\lambda= s^w$ for integer $w$,
then the existence of an $\SOA(\lambda s^t,m,\break s^t,t)$ is
equivalent to that of a $(w,k,m)$-net in base $s$
where $k=w+t$.
\end{lem}

As strong orthogonal arrays are defined without restricting the index
to be a power of $s$,
they provide a more general concept than $(w,k,m)$-nets.
This is in the same spirit as the generalization of orthogonal
Latin squares to orthogonal arrays of strength two.
\citet{HeTan13} discussed several families of strong orthogonal arrays
that cannot be obtained from $(w,k,m)$-nets.
Because of Lemma~\ref{lemma2}, it is not unreasonable to expect
that the star discrepancy of
strong orthogonal arrays would also be $O( n^{-1} (\log n)^{m-1})$
just like nets,
although a precise presentation and rigorous derivation of
this result may require some serious work.
Since our focus is the finite sample space-filling properties of
strong orthogonal arrays,
we choose not to dwell any further on the issue of discrepancy
in this paper.

\section{Characterizing strong orthogonal arrays of strength three}\label{sec3}

Central to our characterizing result is the
notion of embeddability and semi-embeddabili\-ty for orthogonal arrays.

\begin{defin}\label{defin1}
An orthogonal array $\OA(n,m,s,t)$ is said to be embeddable
if it can be obtained by deleting one column from an $\OA(n,m+1,s,t)$.
\end{defin}

Consider the first column of an $\OA(n,m,s,t)$. Then the $s$ levels
in this first column divide the whole array into
$s$ sub-arrays, which are not orthogonal arrays but all become
$\OA(n/s, m-1,s, t-1)$'s if their first columns are deleted.
We say that these $s$ arrays of strength $t-1$
are obtained by branching the first column.
Similarly, branching any other column also produces
$s$ orthogonal arrays of strength $t-1$.
In total, $ms$ such arrays of strength $t-1$ can be obtained.
For easy reference, they are called child arrays or
simply children of the $\OA(n,m,s,t)$ under consideration.

\begin{defin}\label{defin2}
An $\OA(n,m,s,t)$ is said to be semi-embeddable if all of its $ms$
children are
embeddable.
\end{defin}

The following result is immediate.

\begin{lem}\label{lemma3}
If an $\OA(n,m,s,t)$ is embeddable, then it must be semi-embeddable.
\end{lem}

The converse of Lemma~\ref{lemma3} is not always true, and we will see many
examples in the
rest of the paper. One result in \citet{HeTan13} states that
if an embeddable $\OA(n,m,s,3)$ is available, then an $\SOA(n,m,s^3,3)$
can be constructed. The main result of this paper, the following
Theorem~\ref{teo1},
provides a complete characterization for the existence of an $\SOA(n,m,s^3,3)$.

\begin{teo}\label{teo1}
An $\SOA(n,m,s^3,3)$ exists if and only if a semi-embeddable $\OA(n,m,s,3)$
exists.
\end{teo}

The\vspace*{1pt} proof, as given in \hyperref[app]{Appendix}, is actually constructive, and it shows
how to construct an $\SOA(n,m,s^3,3)$ from a semi-embeddable $\OA(n,m,s,3)$
and vice versa. While the characterization is fundamental of
a strong orthogonal array through a generalized orthogonal array
as in \citet{HeTan13},
Theorem~\ref{teo1} does provide a more direct and penetrating characterization
for strong orthogonal arrays of strength three.

As an immediate application of Theorem~\ref{teo1}, we present the following result
on the maximum number of constraints on strong orthogonal arrays.

\begin{teo}\label{teo2}
We have $h(n,s,3) = f(n,s,3)-1$,
provided that
%
\begin{equation}\label{eq3}
f(n,s,3)= f(n/s,s,2)+1,
\end{equation}
where $h(n,s,t)$ and $f(n,s,t)$ are the largest $m$ for an
$\SOA(n,m,s^t,t)$
and an $\OA(n,m,s,t)$ to exist, respectively.
\end{teo}

We know from \citet{HeTan13} that
$ f(n, s, 3) - 1 \leq h(n,s,3) \leq\break  f(n,s,3)$.
Theorem~\ref{teo2} then follows from Theorem~\ref{teo1} if we can show that,
under the condition in (\ref{eq3}),
any $\OA(n,m',s,3)$ with $m' = f(n,s,3)$ is not semi-embeddable. This is
obvious as none of its child arrays, which are $\OA(n/s,m'-1,s,2)$'s,
can be embeddable due to $m'-1 = f(n/s,s,2)$.

For $s=2$, the condition in (\ref{eq3}) is always met, and this
special case of Theorem~\ref{teo2} was obtained in \citet{HeTan13}.
Another important case where the condition in (\ref{eq3})
holds is when $n=s^3$ and $s$ is an even prime
power, in which case we have $f(n,s,3) = s+2$ and $f(n/s,s,2) =
s+1$ [\citet{HedSloStu99}].

The results of \citeauthor{BieEdeSch02} [(\citeyear{BieEdeSch02}), Section~7]
can be regarded as a linear version of Theorem~\ref{teo1}. As such, the following\vadjust{\goodbreak}
Propositions \ref{prop1} and \ref{prop2} have also been established by these authors albeit
in different terminology.

\begin{prop}\label{prop1}
A linear orthogonal array $\OA(s^k, m, s, 3)$ is semi-\break embeddable,
so long as $m \leq(s^{k-1}-1)/(s-1)$.
\end{prop}

An orthogonal array $\OA(n,m,s,t)$ is said to be linear if its runs,
as vectors based on a finite field GF$(s)$, form a linear space.
Proposition~\ref{prop1} can also be established directly. We omit the details
but provide the following\vspace*{1pt} pointers for those readers who are interested in
a direct proof.
Any linear $\OA(s^{k-1}, m-1, s,2)$ is a sub-array
of the saturated linear $\OA(s^{k-1}, (s^{k-1}-1)/(s-1), s,2)$ from
Rao--Hamming construction [\citet{Che14}, Chapter~9].
Permuting\vspace*{1pt} the levels within a column
of this saturated linear array
generates another $\OA(s^{k-1}, (s^{k-1}-1)/(s-1), s,2)$, which is not
linear in general.
Any child of a linear $\OA(s^k, m, s, 3)$ is
an $\OA(s^{k-1}, m-1, s, 2)$, which is either linear or can be obtained
from a linear array by permuting the levels in its columns.

Bush construction gives a linear
$\OA(s^3, s+1, s, 3)$, which can be embedded into an $\OA(s^3, s+2,s,3)$
when $s$ is an even prime power. For odd prime power $s$,
this $\OA(s^3, s+1, s, 3)$ is not embeddable as in this case
$f(s^3,s,3) = s+1$. However, according to Proposition~\ref{prop1}, it is semi-embeddable.
Therefore, an $\SOA(s^3,s+1,s^3,3)$ can always be constructed
when $s$ is a prime power. Examples are
$\SOA(27, 4, 27, 3)$, $\SOA(64, 5, 64, 3)$,
$\SOA(125, 6, 125, 3)$, $\SOA(343, 8, 343, 3)$ and so on.
We summarize the above discussion in the next result.

\begin{prop}\label{prop2}
For any prime power $s$, we have that $h(s^3, s, 3) = s+1$.
\end{prop}

Consider a linear $\OA(s^4, s^2+1, s, 3)$ based on an ovoid; see
Hedayat, Sloane and Stufken [(\citeyear{HedSloStu99}), Section~5.9]. This array satisfies the condition in
Proposition~\ref{prop1} and is therefore semi-embeddable. As such,
an $\SOA(s^4,s^2+1,s^3,3)$ can be constructed by Theorem~\ref{teo1}. Note that
the $\OA(81,10,3,3)$ resulting from taking $s=3$ is not embeddable
as $f(81,3,3)=10$.

If a run occurs more than once in an orthogonal array, it is called
a repeated run. The following Theorem~\ref{teo3} asserts that certain orthogonal arrays
are not semi-embeddable if they have repeated runs.
The proof of Theorem~\ref{teo3} requires the use of a result on orthogonal arrays
with repeated runs, and this is presented in Lemma~\ref{lemma4}.

\begin{lem}\label{lemma4}
If there exists an $\OA(2s^t, m, s,t)$ with a repeated run, then we must
have $m \leq s+t-1$.
\end{lem}

\begin{teo}\label{teo3}
For $s\geq3$, an $\OA(2s^3, s+2, s,3)$ containing
a repeated run is not semi-embeddable.\vadjust{\goodbreak}
\end{teo}

The proofs of Lemma~\ref{lemma4} and Theorem~\ref{teo3} are given in \hyperref[app]{Appendix}.
The bound in Lemma~\ref{lemma4} is quite sharp. For example, taking $t=2$ gives
$m \leq s+1$, which is attainable by the $\OA(2s^2, s+1,s,2)$
from juxtaposing
two identical $\OA(s^2,s+1,s,2)$'s where $s$ is a prime power.

\section{Construction of $\operatorname{SOA}(54,5,27,3)$}\label{sec4}

In the present section, we discuss the application of the results in
Section~\ref{sec3}
to the existence and construction of $\SOA(54,5,27,3)$'s.
According to \citet{HedSeiStu97}, the maximum number
$m$ of factors in an orthogonal array $\OA(54, m, 3, 3)$ is five and
there are exactly four nonisomorphic $\OA(54, 5, 3, 3)$'s.
These four arrays, labeled as I, II, III and IV, are available
in explicit form in their paper.
To study the existence and construction of $\SOA(54,5,27,3)$'s, Theorem~\ref{teo1}
says that it suffices to examine the semi-embeddability of these
four nonisomorphic $\OA(54, 5, 3, 3)$'s.

Array I has two repeated runs and array II has one repeated run.
By Theorem~\ref{teo3}, neither array is semi-embeddable. Thus, no $\SOA(54, 5, 27,3)$
can be constructed from array I or array II.
Both arrays III and IV have no repeated run. It is thus possible for them
to be semi-embeddable. Our direct computer search shows that this is
indeed the case. The two $\SOA(54, 5, 27, 3)$'s constructed from
these two arrays using Theorem~\ref{teo1} are given in Tables~\ref{table1} and \ref{table2},
respectively. To save space,
both of them are presented in transposed forms,
with the top half of each table displaying runs 1--27 and the bottom half
runs 28--54.

\begin{sidewaystable}
\tablewidth=\textwidth
\tabcolsep=0pt
\caption{$\SOA(54, 5, 27,3)$ constructed from the third $\OA(54, 5, 3, 3)$, array III}\label{table1}
\begin{tabular*}{\tablewidth}{@{\extracolsep{\fill}}@{}ld{2.0}d{2.0}d{2.0}d{2.0}d{2.0}d{2.0}d{2.0}d{2.0}d{2.0}d{2.0}d{2.0}d{2.0}
d{2.0}d{2.0}d{2.0}d{2.0}d{2.0}d{2.0}d{2.0}d{2.0}d{2.0}d{2.0}d{2.0}d{2.0}d{2.0}d{2.0}@{}}
\hline
\phantom{0}0 & 3 & 0 & 8 & 8 & 18 & 24 & 16 & 3 & 7 & 7 & 9 & 12 & 23 & 6 & 6 &1 & 5 & 1 & 5 & 9 & 18 & 15 & 21 & 13 & 26 & 4\\
\phantom{0}0 & 3 & 1 & 18 & 26 & 6 & 8 & 15 & 5 & 9 & 13 & 6 & 7 & 21 & 7 & 8 &9 & 18 & 17 & 22 & 0 & 3 & 2 & 4 & 12 & 24 & 13\\
\phantom{0}0 & 3 & 10 & 8 & 21 & 5 & 24 & 6 & 20 & 1 & 9 & 7 & 9 & 6 & 13 & 23& 8 & 4 & 18 & 15 & 2 & 7 & 24 & 15 & 3 & 0 & 17\\
\phantom{0}1 & 5 & 10 & 22 & 7 & 25 & 4 & 7 & 20 & 11 & 2 & 11 & 8 & 8 & 12 &21 & 18 & 15 & 6 & 3 & 24 & 15 & 0 & 6 & 3 & 0 & 25\\
\phantom{0}9 & 18 & 9 & 15 & 15 & 14 & 14 & 10 & 18 & 21 & 21 & 25 & 25 & 20 &0 & 0 & 3 & 6 & 3 & 6 & 7 & 5 & 7 & 5 & 1 & 2 & 12
\\[6pt]
\phantom{0}4 & 12 & 15 & 14 & 22 & 11 & 19 & 2 & 2 & 24 & 21 & 25 & 17 & 19 &11 & 10 & 23 & 26 & 20 & 16 & 13 & 10 & 17 & 14 & 20 & 22 & 25\\
 17 & 4 & 5 & 21 & 12 & 19 & 10 & 26 & 22 & 2 & 1 & 15 & 24 & 11 & 20& 11 & 23 & 25 & 19 & 16 & 14 & 10 & 25 & 23 & 20 & 14 & 16\\
 25 & 17 & 19 & 1 & 4 & 12 & 12 & 25 & 14 & 22 & 14 & 2 & 5 & 21 & 18& 23 & 20 & 10 & 16 & 16 & 26 & 13 & 11 & 22 & 26 & 19 & 11\\
 16 & 19 & 16 & 13 & 13 & 1 & 4 & 14 & 26 & 14 & 23 & 23 & 20 & 2 & 5& 22 & 19 & 10 & 17 & 17 & 26 & 12 & 21 & 9 & 24 & 9 & 18\\
 12 & 16 & 16 & 13 & 17 & 13 & 17 & 24 & 24 & 23 & 23 & 26 & 22 & 26& 22 & 10 & 11 & 11 & 20 & 19 & 19 & 1 & 4 & 4 & 2 & 8 & 8\\
\hline
\end{tabular*}
\vspace*{20pt}
\tablewidth=\textwidth
\tabcolsep=0pt
\caption{$\SOA(54, 5, 27,3)$ constructed from the fourth $\OA(54, 5, 3, 3)$, array IV}\label{table2}
\begin{tabular*}{\tablewidth}{@{\extracolsep{\fill}}@{}ld{2.0}d{2.0}d{2.0}d{2.0}d{2.0}d{2.0}d{2.0}d{2.0}d{2.0}d{2.0}d{2.0}d{2.0}
d{2.0}d{2.0}d{2.0}d{2.0}d{2.0}d{2.0}d{2.0}d{2.0}d{2.0}d{2.0}d{2.0}d{2.0}d{2.0}d{2.0}@{}}
\hline
\phantom{0}0 & 3 & 0 & 8 & 8 & 18 & 24 & 16 & 3 & 7 & 7 & 9 & 12 & 23 & 6 & 6 &1 & 5 & 1 & 5 & 9 & 18 & 15 & 21 & 13 & 26 & 4\\
\phantom{0}0 & 3 & 1 & 18 & 26 & 6 & 8 & 15 & 5 & 9 & 13 & 6 & 7 & 21 & 7 & 8 &9 & 18 & 17 & 22 & 0 & 3 & 2 & 4 & 12 & 24 & 16\\
\phantom{0}0 & 3 & 10 & 8 & 21 & 5 & 24 & 6 & 23 & 1 & 9 & 7 & 9 & 6 & 14 & 19 &8 & 4 & 18 & 15 & 2 & 7 & 24 & 12 & 3 & 0 & 13\\
\phantom{0}1 & 5 & 10 & 22 & 7 & 25 & 4 & 7 & 23 & 11 & 2 & 11 & 8 & 8 & 18 & 12& 18 & 15 & 6 & 3 & 24 & 12 & 0 & 6 & 3 & 0 & 13\\
\phantom{0}9 & 18 & 9 & 15 & 15 & 14 & 14 & 10 & 18 & 21 & 21 & 25 & 25 & 20 & 0& 0 & 3 & 6 & 3 & 6 & 7 & 5 & 7 & 5 & 1 & 2 & 12
\\[6pt]
\phantom{0}4 & 12 & 15 & 14 & 22 & 11 & 19 & 2 & 2 & 21 & 24 & 25 & 17 & 19 & 11& 10 & 23 & 26 & 20 & 16 & 13 & 10 & 17 & 14 & 25 & 20 & 22\\
14 & 4 & 5 & 21 & 12 & 19 & 10 & 26 & 22 & 2 & 1 & 15 & 24 & 11 & 20& 11 & 23 & 25 & 19 & 13 & 17 & 10 & 25 & 23 & 14 & 20 & 16\\
26 & 17 & 22 & 1 & 4 & 12 & 15 & 25 & 17 & 20 & 16 & 2 & 5 & 21 & 18& 20 & 19 & 11 & 13 & 14 & 25 & 16 & 10 & 23 & 22 & 26 & 11\\
25 & 22 & 16 & 13 & 16 & 1 & 4 & 17 & 26 & 20 & 17 & 23 & 20 & 2 & 5& 19 & 10 & 19 & 14 & 26 & 14 & 15 & 9 & 21 & 9 & 24 & 21\\
12 & 16 & 16 & 13 & 17 & 13 & 17 & 24 & 24 & 23 & 23 & 26 & 22 & 26 &22 & 10 & 11 & 11 & 20 & 19 & 19 & 1 & 4 & 4 & 8 & 2 & 8\\
\hline
\end{tabular*}
\end{sidewaystable}

From Proposition~\ref{prop2}, we know that $h(s^3,s,3) =s+1$
if $s$ is a prime
power. Construction of $\SOA(54, 5, 27, 3)$
establishes the following result.

\begin{teo}\label{teo4}
We have that $h(54, 3, 3) = f(54, 3, 3) = 5$.
\end{teo}

To gain some insights into the semi-embeddability and
nonsemi-embedd\-ability of the four
nonisomorphic $\OA(54, 5, 3, 3)$'s, we make use of the enumeration
\mbox{results} on orthogonal arrays of 18 runs.
\citet{Sch09} enumerated all nonisomorphic orthogonal arrays of 18
runs and
he found that there are exactly 12 nonisomorphic $\OA(18,4,3,2)$'s,
which he labeled as $4.0.i$ for $i =1, \ldots, 12$ in his paper.
Among these 12 arrays, five of them are nonembeddable and the other seven
are embeddable. The nonembeddable ones are 4.0.5, 4.0.7, 4.0.10, 4.0.11
and 4.0.12.

For a given $\OA(54, 5, 3, 3)$, three child arrays can be obtained by
branching each of the five columns.
For the first $\OA(54, 5, 3, 3)$, array I, among the three child arrays
from branching each column, one is isomorphic to 4.0.1 and the other two
are isomorphic to 4.0.5. For array II,
one of the three arrays from branching each column is isomorphic to
4.0.5 and the other two are isomorphic to 4.0.2.
Thus, the early conclusion that arrays I and II are not semi-embeddable
can also be drawn from the fact that array 4.0.5 is not embeddable.
For array III,
all the 12 child arrays from branching columns 1 through 4 are isomorphic
to 4.0.4, and the three child arrays from branching column 5 are isomorphic
to 4.0.1. For array IV, the 9 child arrays from branching columns 1, 2 and
5 are isomorphic to 4.0.2, and the 6 child arrays from branching columns
3 and 4 are isomorphic to 4.0.4. Since all these child arrays are embeddable,
arrays III and IV are semi-embeddable.

\section{Discussion and future work}\label{sec5}

He and Tang [(\citeyear{HeTan13}), Theorem~1] presented a general method of constructing
strong orthogonal arrays from ordinary orthogonal arrays.
For the case of strength three, this result means that an $\SOA(n,m, s^3,3)$
can be constructed from an $\OA(n,m+1,s,3)$.
Specifically, let $(a_1, \ldots, a_m, a_{m+1})$ be an $\OA(n,m+1,s,3)$.
Then $B=\{(a_1, b_1, c_1);\break  \ldots; (a_m,b_m,c_m)\}$ is a  $\GOA(n,m,s,3)$
and $D=(d_1, \ldots, d_m)$ is an $\SOA(n,m,\break s^3,3)$,
where $(b_1, \ldots, b_m) = (a_{m+1}, \ldots, a_{m+1})$,
$(c_1, \ldots, c_m) = (a_2, \ldots, a_m, a_1)$, and
$d_i = a_i s^2 + b_i s+ c_i$ for $i=1,\ldots,m$.
Note that the same $a_{m+1}$ is taken for all $b_i$'s.
This simple construction of strong orthogonal arrays
from orthogonal arrays should be sufficient
for most practical purposes because only one column is lost during the
construction. In the terminology of the present paper,
semi-embeddability of $A=(a_1, \ldots, a_m)$ is automatic due to
its embeddability into $(a_1, \ldots, a_m, a_{m+1})$.

The above construction also says that embeddability of an $\OA(n,m,s,3)$
is sufficient for the existence of an $\SOA(n,m,s^3,3)$. The present paper
strengthens this result by proving that an $\SOA(n,m,s^3,3)$ exists
if and only if a semi-embeddable $\OA(n,m,s,3)$ exists. The path
of constructing an $\SOA(n,m,s^3,3)$ from a semi-embeddable $\OA(n,m,s,3)$
is still via a $\GOA(n,m,s,3)$, but
there is some difference. To illustrate, let $A=(a_1, \ldots, a_m)$
be a semi-embeddable $\OA(n,m,s,3)$. Then in forming
$B=\{(a_1, b_1, c_1); \ldots;\break  (a_m,b_m,c_m)\}$,
a $\GOA(n,m,s,3)$,
although we still take
$(c_1, \ldots, c_m) = (a_2, \ldots,\break  a_m, a_1)$,
the columns $b_i$'s as obtained in the proof of Theorem~\ref{teo1} in the \hyperref[app]{Appendix}
cannot be all the same unless the semi-embeddable $A=(a_1, \ldots, a_m)$
is also embeddable. This is because
if the $b_i$'s equal the same column, say $b$, then
$(a_1, \ldots, a_m, b)$ must be an $\OA(n,m+1,s,3)$ due to the fact
that $B=\{(a_1, b, c_1); \ldots; (a_m,b,c_m)\}$ is
a $\GOA(n,m,s,3)$.

The question arises of if a given orthogonal array is semi-embeddable.
The simplest case is to consider sub-arrays of
the available orthogonal arrays by deleting one or more columns.
All arrays obtained this way are embeddable, and hence semi-embeddable.
In Section~\ref{sec3}, we have presented two further results for judging
whether or not an orthogonal array is semi-embeddable.
Proposition~\ref{prop1} tells us that a linear orthogonal array $\OA(s^k, m, s,
3)$ is
semi-embeddable, provided $m \leq(s^{k-1}-1)/(s-1)$. This result
has led to the conclusion that
$\OA(s^3, s+1, s, 3)$ from Bush construction and
$\OA(s^4, s^2+1, s,3)$ base on an ovoid are both
semi-embeddable, where $s$ is any prime power.
Theorem~\ref{teo3} states that an $\OA(2s^3, s+2,s,3)$ is
not semi-embeddable for $s \geq3$ if it has a repeated run, allowing
us to immediately identify two nonsemi-embeddable $\OA(54, 5, 3, 3)$'s
in Section~\ref{sec4}. When none of the above methods can give a definitive answer,
one can make use of relevant enumeration results if they are available or
conduct a complete search as a last resort.
In Section~\ref{sec4}, we have done it both ways in
determining the semi-embeddability of the other two $\OA(54, 5, 3, 3)$'s.

One obvious future direction is to study to what extent the current work
can be extended to strong orthogonal arrays of strength four or higher.
Although such extension work may not be as neat as what we have done
for strong orthogonal arrays of strength three, some
useful results are still possible. We leave this to the future.

A more promising direction is what can be done when orthogonal arrays of
strength three or higher are too expensive to use for given resources.
As discussed in \citet{HeTan13}, strong orthogonal arrays of strength
two can be straightforwardly constructed from ordinary orthogonal arrays
of strength two but the former do not improve upon the latter in terms of
lower dimensional space-filling. The question then is if we can construct
designs that, although not strong orthogonal arrays of strength three,
are better than strong orthogonal arrays
of strength two. Some preliminary results have been obtained,
and we hope to write a future paper along this direction.

\begin{appendix}\label{app}
\section*{Appendix}

\begin{pf*}{Proof of Theorem~\ref{teo1}}
We first prove that the existence of an $\SOA(n,\break m,s^3,3)$ implies
the existence of a semi-embeddable $\OA(n,m,s,3)$.
Suppose there exists an $\SOA(n,m,s^3,3)$. Then Lemma~\ref{lemma1} implies the
existence of
a $\GOA(n,m,s,3)$, $B=\{(a_1,b_1,c_1); \ldots; (a_m,b_m,c_m) \}$.
We will show
that array $A=(a_1, \ldots, a_m)$ is a semi-embeddable $\OA(n,m,s,3)$.
That $A$ is an $\OA(n,m,s,3)$ follows directly from the definition of
generalized orthogonal arrays. By Definition~\ref{defin2}, what remains to be shown
is that the children of $A$ are all embeddable. Now consider array
$P = (a_1, \ldots, a_m, b_1)$. That $B=\{(a_1,b_1,c_1);
\ldots; (a_m,b_m,c_m) \}$ is a $GOA(n,m,s,3)$ dictates that
$(a_1, a_j, b_1)$ is an $\OA(n,3,s,3)$ for any $j=2, \ldots, m$. This
implies that the array $Q$ obtained by selecting the $n/s$ rows of
$(a_2, \ldots, a_m, b_1)$ that correspond to a given level in $a_1$
must be an $\OA(n/s, m, s, 2)$. Clearly, array $Q$ becomes a child of
$A$ if
the last column is deleted. This shows that all the $s$ children
of $A$ from branching column $a_1$ are embeddable. The same argument
also applies to the children from branching other columns of $A$.

We next show that an $\SOA(n,m,s^3,3)$ can be constructed
from a semi-embeddable $\OA(n,m,s,3)$.
Suppose that $A=(a_1, \ldots, a_m)$ be a semi-embedd\-able $\OA(n,m,s,3)$.
We will construct a $\GOA(n,m,s,3)$,
$B=\{(a_1,b_1,c_1); \ldots;\break  (a_m,b_m,c_m) \}$.
Then Lemma~\ref{lemma1} allows an $\SOA(n,m,s^3,3)$ to be constructed from~$B$.
The last paragraph shows that if $(a_1, a_j, b_1)$ for any $j=2, \ldots,m$
is an $\OA(n,3,s,3)$, then all the children of $A$ from branching column
$a_1$ are embeddable. We observe that this argument is
entirely reversible, meaning that if all the children of $A$
from branching column $a_1$ are embeddable, then a column $b_1$ can be
obtained so that $(a_1, a_j, b_1)$ is an $\OA(n,3,s,3)$
for any $j=2, \ldots, m$. Similarly, a column $b_i$ for $i=2,\ldots, m$
can be obtained so that $(a_i, a_j, b_i)$ is an $\OA(n,3,s,3)$ for
any $j=1,\ldots, i-1, i+1, \ldots,m$.
Take $(c_1, \ldots, c_m) = (a_2,\ldots, a_m, a_1)$. Now it is evident
that array $B=\{(a_1,b_1,c_1); \ldots; (a_m,b_m,c_m) \}$ is a $\GOA(n,m,s,3)$.
\end{pf*}

\begin{pf*}{Proof of Lemma~\ref{lemma4}}
Let $c$ be a repeated run of an $\OA( 2 s^t, m, s, t)$.
For $i=0,1,\ldots, m$, let
$n_i$ be the number of other runs that have exactly $i$ coincidences
with $c$. As $c$ is a repeated run, we must have $n_m \geq1$.
A result from \citet{BosBus52} states that
\begin{equation}\label{A1}
\sum_{i=j}^m \comb{i} {j}
n_i = \comb{m} {j} \bigl( 2 s^{t-j}-1 \bigr)\qquad\mbox{where }j=0, 1, \ldots, t.
\end{equation}
Choosing $j=t$ in (\ref{A1}) gives $\sum_{i=t}^m \tcomb{i}{t} n_i = \tcomb{m}{t}$.
Combining this equation with the fact that $n_m \geq1$, we must have
$n_t = \cdots= n_{m-1} =0$ and $n_m =1$.
Now consider the two equations given by
setting $j=t-1$ and $j=t-2$ in (\ref{A1}).
Solving these two equations, we obtain
$n_{t-1} = 2 (s-1) \tcomb{m}{t-1}$ and $n_{t-2} = 2(s-1) (s+t-1-m)
\tcomb{m}{t-2}$. As $n_{t-2} \geq0$, we must have
$s+t-1-m \geq0$, implying that $m \leq s+t-1$. Lemma~\ref{lemma4} is proved.
\end{pf*}

\begin{pf*}{Proof of Theorem~\ref{teo3}}
An $\OA(2s^3,s+2,s,3)$ containing a repeated run gives rise to many
children that have repeated runs. Let $A_0$ be any such child array
$\OA(2s^2, s+1, s,2)$ with a repeated run. By Definition~\ref{defin2},
Theorem~\ref{teo3} will be established if we can show that $A_0$ is not embeddable.
Let $c$ be a repeated run
of $A_0$ and $n_i$ be the number of other runs of $A_0$ that have
exactly $i$ coincidences with~$c$. Applying the results in the proof
of Lemma~\ref{lemma4} to the case $t=2$ and $m=s+1$, we obtain
\begin{equation}\label{A2}
n_0 = n_2 = \cdots= n_s =0,\qquad
n_1 = 2 \bigl(s^2-1 \bigr),\qquad n_{s+1}= 1.
\end{equation}

We will prove that $A_0$ is not embeddable by the method of contradiction.
Suppose that $A_0$ is embeddable and let $A_0^+$ be an $\OA(2s^2, s+2,s,2)$
obtained from $A_0$ by adding one column. Recall that $c$ is a repeated
run of $A_0$. Let $c^+$ be the run of $A_0^+$ corresponding to $c$.
Let $n_i^+$ be the number of other runs of $A_0^+$ that have exactly
$i$ coincidences with $c^+$. By Lemma~\ref{lemma4}, array $A_0^+$ cannot have
a repeated run, implying that $n_{s+2}^+ =0$. Noting that
$A_0$ is a sub-array of $A_0^+$, and combining $n_{s+2}^+ =0$ with
the results in (\ref{A2}), we obtain
$n_0^+ = n_3^+ = \cdots= n_s^+ =0$, $n_{s+1}^+ =1$
and
\begin{equation}\label{A3}
n_1^+ + n_2^+ = n_1 = 2
\bigl(s^2-1 \bigr).
\end{equation}
On the other hand, the coincidence equation in (\ref{A1}) becomes
\begin{equation}\label{A4}
\sum_{i=j}^{s+2} \comb{i} {j}
n_i^+ = \comb{s+2} {j} \bigl( 2 s^{2-j}-1 \bigr)\qquad\mbox{where }j=0, 1, 2.
\end{equation}
Using the two equations from taking $j=1,2$ in (\ref{A4}) and
the already obtained results about $n_i^+$ for $i =3, \ldots, s+2$,
we obtain\vspace*{1pt}
$n_1^+ = 2 s^2 -5 $ and $n_2^+ = s+1$, which gives
$n_1^+ + n_2^+ = 2 s^2 + s - 4$. But this contradicts
(\ref{A3}) for $s \geq3$ because
$ (2s^2+s -4) - 2(s^2-1) = s-2 \geq1$ for any $s \geq3$.
The proof is complete.
\end{pf*}
\end{appendix}




\printaddresses
\end{document}